\documentclass[12pt]{article}
= 0.0in \topmargin = 0.3in \oddsidemargin=0.1in
\def\Dj{\hbox{D\kern-.73em\raise.30ex\hbox{-}
\raise-.30ex\hbox{}}}
\def\dj{\hbox{d\kern-.33em\raise.80ex\hbox{-}
\raise-.80ex\hbox{\kern-.40em}}}
\usepackage{epsfig}
\usepackage{amsmath,amsthm,amsfonts,amssymb,amscd}
\usepackage{graphicx}

\newtheorem{theorem}{Theorem}[section]
 \newtheorem{example}{Example}
 \newtheorem{problem}[example]{Problem}
 \newtheorem{defin}[theorem]{Definition}

 \newtheorem{nt}{Note}
 \newtheorem{conj}{Conjecture}

 \newcommand{\singlespacing}{\let\CS=\@currsize\renewcommand{\baselinestretch}{1}\tiny\CS}
 \newcommand{\oneandahalfspacing}{\let\CS=\@currsize\renewcommand{\baselinestretch}{1.25}\tiny\CS}
 \newcommand{\doublespacing}{\let\CS=\@currsize\renewcommand{\baselinestretch}{1.35}\tiny\CS}

 \newtheorem{rule-def}[theorem]{Rule}

\begin{document}

\baselineskip=0.30in

\vspace*{20mm}

\begin{center}
{\Large \bf Wiener Indices of Spiro and Polyphenyl Hexagonal Chains
$^\dagger$}

\vspace{10mm}

{\large \bf Hanyuan Deng$^\ddagger$}

\vspace{6mm}

\baselineskip=0.20in

{\it College of Mathematics and Computer Science,
\\ Hunan Normal University, Changsha, Hunan 410081, P. R.
China\/}

\vspace{6mm}


\end{center}

\vspace{6mm}

\baselineskip=0.20in

\noindent {\bf Abstract }

\vspace{3mm}

{\small The Wiener index $W(G)$ of a connected graph $G$ is the sum
of distances between all pairs of vertices in $G$. In this paper, we
first give the recurrences or explicit formulae for computing the
Wiener indices of spiro and polyphenyl hexagonal chains, which are
graphs of a class of unbranched multispiro molecules and polycyclic
aromatic hydrocarbons, then we establish a relation between the
Wiener indices of a spiro hexagonal chain and its corresponding
polyphenyl hexagonal chain, and determine the extremal values and
characterize the extremal graphs with respect to the Wiener index
among all spiro and polyphenyl hexagonal chains with $n$ hexagons,
respectively. An interesting result shows that the average value of
the Wiener indices with respect to the set of all such hexagonal
chains is exactly the average value of the Wiener indices of three
special hexagonal chains, and is just the Wiener index of the
meta-chain.}

\vspace{30mm}

\baselineskip=0.20in

\vspace{5mm}

\noindent
-------------------------------------- \\
$^\dagger${\small Project supported by the Scientific Research Fund
of Hunan Provincial Education Department (09A057) and the Hunan
Provincial Natural Science Foundation of China (09JJ6009).}
\\[2mm]
$^\ddagger${\small 
e-mail: {\tt hydeng@hunnu.edu.cn}\,.}

\newpage

\baselineskip=0.30in

\section{Introduction}

All graphs considered in this paper are simple, undirected and
connected. The vertex and edge sets of a graph $G$ are $V(G)$ and
$E(G)$, respectively. The distance $d_G(u, v)$ between vertices $u$
and $v$ is the number of edges on a shortest path connecting these
vertices in $G$. The distance $W(G,v)$ of a vertex $v\in V(G)$ is
the sum of distances between $v$ and all other vertices of $G$.

The Wiener index [1,2] of a graph $G$ is a graph invariant based on
distances. It is defined as the sum of distances between all pairs
of vertices in $G$:
$$W(G)=\sum\limits_{\{u,v\}\subseteq V(G)}d_G(u, v)=\frac{1}{2}\sum\limits_{v\in V(G)}W(G, v).$$

Wiener index is the oldest topological index related to molecular
branching [3]. A tentative explanation of the relevance of the
Wiener index in research of QSPR and QSAR is that it correlates with
the van der Waals surface area of the molecule [4]. Until now, the
Wiener index has gained much popularity and new results related to
it are constantly being reported. For a survey of results and
further bibliography on the chemical applications and the
mathematical literature of the Wiener index, see [5-8] and the
references cited therein.

Spiro compounds are an important class of cycloalkanes in organic
chemistry. A spiro union in spiro compounds is a linkage between two
rings that consists of a single atom common to both rings and a free
spiro union is a linkage that consists of the only direct union
between the rings. The common atom is designated as the spiro atom.
According to the number of spiro atoms present, compounds are
distinguished as monospiro, dispiro, trispiro, etc, ring systems.
Figure 1(i) illustrates three linear polyspiro alicyclic
hydrocarbons. Here, we consider a subclass of unbranched multispiro
molecules, in which every ring is a hexagon, and their graphs are
called spiro hexagonal chains (or chain hexagonal cacti, or
six-membered ring spiro chain [9-11]).

Two or more benzene rings are linked by a cut edges consisting of
aromatics called polycyclic aromatic hydrocarbons which is a class
of aromatics. A class of compounds in which two and more benzene
rings are directly linked by a cut edge known as the biphenyl
compounds, and their graphs are called polyphenyl hexagonal chains
[12]. Figure 1(ii) illustrates ortho-terphenyl, meta-terphenyl and
pera-terphenyl.

Some explicit recurrences for the matching and independence
polynomials in the spiro and polyphenyl hexagonal chains were
derived in [9], and the spiro and polyphenyl hexagonal chains with
the extremal values of the Merrifield-Simmons index and Hosoya index
were determined in [11,12]. The extremal energies of the spiro and
polyphenyl hexagonal chains were found in [10].

In this paper, we will first give the recurrences or explicit
formulae for computing the Wiener indices of spiro and polyphenyl
hexagonal chains, and then establish a relation between a spiro
hexagonal chain and its corresponding polyphenyl hexagonal chain and
determine the extremal values and characterize the extremal graphs
with respect to the Wiener index among all spiro hexagonal chains
and polyphenyl hexagonal chains with $n$ hexagons. Also, we will
discuss the average value of the Wiener indices with respect to the
set of all such hexagonal chains, and find an interesting result
which shows that the average value is exactly the average value of
three special hexagonal chains, and is just the Wiener index of the
meta-chain.

\section{Wiener index of spiro hexagonal chains}

A hexagonal cactus is a connected graph in which every block is a
hexagon. A vertex shared by two or more hexagons is called a
cut-vertex. If each hexagon of a hexagonal cactus $G$ has at most
two cut-vertices, and each cut-vertex is shared by exactly two
hexagons, then $G$ is called a spiro hexagonal chain. The number of
hexagons in $G$ is called the length of $G$. An example of a spiro
hexagonal chain is shown in Figure 2(i).

Obviously, a spiro hexagonal chain of length $n$ has $5n+1$ vertices
and $6n$ edges. Furthermore, any spiro hexagonal chain of length
greater than one has exactly two hexagons with only one cut-vertex.
Such hexagons are called terminal hexagons. Any remaining hexagons
are called internal hexagons.

Let $G_n=H_0H_1\cdots H_{n-1}$ be a spiro hexagonal chain of length
$n (n\geq 3)$. $H_k$ is the $(k+1)$-th hexagon of $G_n$, $c_k$ is
the common cut-vertex of $H_{k-1}$ and $H_k$, $k=1,2,\cdots,n-1$.
Then, the sequence $(c_2, c_3, \cdots, c_{n-1})$ of length $n-2$ is
called the cut-vertex sequence of $G_n$. Obviously, $G_n$ is
determined completely by its cut-vertex sequence. A vertex $v$ of
$H_k$ is said to be ortho-, meta- and para-vertex of $H_k$ if the
distance between $v$ and $c_k$ is 1, 2 and 3, denoted by $o_k$,
$m_k$ and $p_k$, respectively. Examples of ortho-, meta-, and
para-vertices are shown in Figure 3(ii). Except the first hexagon,
any hexagon in a spiro hexagonal chain has two ortho-vertices, two
meta-vertices and one para-vertex.

A spiro hexagonal chain $G_n$ is a spiro ortho-chain if
$c_k=o_{k-1}$ for $2\leq k\leq n-1$, i.e., its cut-vertex sequence
is $(o_1, o_2, \cdots, o_{n-2})$. The spiro meta-chain and spiro
para-chain are defined in a completely analogous manner. The spiro
ortho-, meta- and para-chain of length $n$ is denoted by $O_n$,
$M_n$ and $P_n$, respectively. Examples of spiro ortho-, meta-, and
para-chains are shown in Figure 4.

In the following, we first give a recurrence for computing the
Wiener indices of spiro hexagonal chains, and then derive a formula
for computing the Wiener indices of spiro hexagonal chains.

Let $G_n=H_0H_1\cdots H_{n-1}$ be a spiro hexagonal chain with $n$
hexagons as shown in Figure 3(i). $G_{n-1}=H_0H_1\cdots H_{n-2}$ and
$c_{n-1}, o_{n-1}, m_{n-1}, p_{n-1}$ are the cut-, ortho-, meta-,
and para-vertex in $H_{n-1}$, respectively. Then

$$
\begin{array}{ll}
W(G_n, o_{n-1})&=\sum\limits_{v\in G_n}d(v,
o_{n-1})=\sum\limits_{v\in G_{n-1}}(d(v, c_{n-1})+1)+8 \\
&=W(G_{n-1}, c_{n-1})+5(n-1)+9;\\
\end{array}
$$

$$
\begin{array}{ll}
W(G_n, m_{n-1})&=\sum\limits_{v\in G_n}d(v,
m_{n-1})=\sum\limits_{v\in G_{n-1}}(d(v, c_{n-1})+2)+7\\
&=W(G_{n-1}, c_{n-1})+10(n-1)+9;\\
\end{array}
$$

$$
\begin{array}{ll}
W(G_n, p_{n-1})&=\sum\limits_{v\in G_n}d(v,
p_{n-1})=\sum\limits_{v\in G_{n-1}}(d(v, c_{n-1})+3)+6\\
&=W(G_{n-1}, c_{n-1})+15(n-1)+9.\\
\end{array}
$$

So, we have
\begin{equation}
W(G_n, c_n)=W(G_{n-1},c_{n-1})+f(c_n)
\end{equation}
where
$$
f(c_n)=\left\{
\begin{array}{ll}
5(n-1)+9, & c_n \,\,\, is \,\,\, o_{n-1};\\
10(n-1)+9, & c_n \,\,\, is \,\,\, m_{n-1};\\
15(n-1)+9, & c_n \,\,\, is \,\,\, p_{n-1}.\\
\end{array}
\right.
$$

\vspace{3mm}

Let $v_1, v_2, v_3, v_4, v_5$ be the vertices of $H_{n-1}$ different
from $c_{n-1}$. By the definition of Wiener index,
$$
\begin{array}{ll}
W(G_n)&=W(G_{n-1})+\sum\limits_{i=1}^5\sum\limits_{v\in G_{n-1}}d(v,
v_i)+\sum\limits_{1\leq i<j\leq 5}d(v_i, v_j)\\
&=W(G_{n-1})+\sum\limits_{v\in G_{n-1}}(5d(v,
c_{n-1})+9)+18\\
&=W(G_{n-1})+5W(G_{n-1}, c_{n-1})+45n-18.\\
\end{array}
$$
i.e.,
\begin{equation}
W(G_n)=W(G_{n-1})+5W(G_{n-1}, c_{n-1})+45n-18.
\end{equation}

\vspace{3mm}

Combining equations (1) and (2), we can get the following recurrence
for computing the Wiener indices of spiro hexagonal chains.

{\bf Theorem 2.1}. Let $G_n=H_0H_1\cdots H_{n-1}$ be a spiro
hexagonal chain of length $n$, $c_{n-1}$ the cut-vertex of
$H_{n-1}$. Then
$$
\left\{
\begin{array}{l}
W(G_n)=W(G_{n-1})+5W(G_{n-1}, c_{n-1})+45n-18, \\
W(G_n, c_n)=W(G_{n-1},c_{n-1})+f(c_n),\\
W(G_1)=27, \\
W(G_1, c_1)=f(c_1)=9,\\
\end{array}
\right.
$$ and
$$
f(c_n)=\left\{
\begin{array}{ll}
5(n-1)+9, & c_n \,\,\, is \,\,\, o_{n-1};\\
10(n-1)+9, & c_n \,\,\, is \,\,\, m_{n-1};\\
15(n-1)+9, & c_n \,\,\, is \,\,\, p_{n-1}.\\
\end{array}
\right.
$$

\vspace{3mm}

Using the recurrence for computing the Wiener indices of spiro
hexagonal chains, we have
$$
\begin{array}{ll}
W(G_n)&=W(G_{n-1})+5W(G_{n-1}, c_{n-1})+45n-18\\
&=W(G_1)+5\sum\limits_{k=1}^{n-1}W(G_{n-k},
c_{n-k})+45\sum\limits_{k=1}^{n-1}(n-k+1)-18(n-1)\\
&=27+5\sum\limits_{k=1}^{n-1}W(G_{k},
c_{k})+\frac{45}{2}(n+2)(n-1)-18(n-1)\\
&=5\sum\limits_{k=1}^{n-1}W(G_{k}, c_{k})+\frac{45}{2}n^2+\frac{9}{2}n\\
\end{array}
$$
i.e.,
\begin{equation}
W(G_n)=5\sum\limits_{k=1}^{n-1}W(G_{k},
c_{k})+\frac{45}{2}n^2+\frac{9}{2}n.
\end{equation}

\vspace{3mm}

And,

$$
\begin{array}{ll}
W(G_k, c_k)&=W(G_{k-1}, c_{k-1})+f(c_k)\\
&=W(G_1, c_1)+f(c_2)+\cdots+f(c_k)\\
&=f(c_1)+f(c_2)+\cdots+f(c_k),
\end{array}
$$
i.e.,
\begin{equation}
W(G_k, c_k)=\sum\limits_{i=1}^{k}f(c_{i}).
\end{equation}

Combining equations (3) and (4), we can obtain the following formula
for computing the Wiener indices of spiro hexagonal chains.

\vspace{3mm}

{\bf Theorem 2.2}. Let $G_n=H_0H_1\cdots H_{n-1}$ be a spiro
hexagonal chain with $n$ hexagons. $c_{k}$ is the common cut-vertex
of $H_{k-1}$ and $H_{k}$, $1\leq k\leq n-1$. Then
$$
\begin{array}{ll}
W(G_n)&=5\sum\limits_{k=1}^{n-1}\sum\limits_{i=1}^kf(c_{i})+\frac{45}{2}n^2+\frac{9}{2}n\\
&=5\sum\limits_{k=1}^{n-1}(n-k)f(c_k)+\frac{45}{2}n^2+\frac{9}{2}n\\
\end{array}
$$
where
\begin{equation}
f(c_k)=\left\{
\begin{array}{ll}
5(k-1)+9, & c_k\,\,\, is \,\,\, o_{k-1};\\
10(k-1)+9, & c_k \,\,\, is \,\,\, m_{k-1};\\
15(k-1)+9, & c_k \,\,\, is \,\,\, p_{k-1}.\\
\end{array}
\right.
\end{equation}
and $f(c_1)=9$.

\vspace{3mm}

In the spiro orth-chain $O_n$, the spiro meta-chain $M_n$ and the
spiro para-chain $P_n$, $c_k$ is $o_{k-1}, m_{k-1}$ and $p_{k-1}$,
respectively. Then

$W(O_n)=5\sum\limits_{k=1}^{n-1}(n-k)(5(k-1)+9)+\frac{45}{2}n^2+\frac{9}{2}n=\frac{25}{6}n^3+\frac{65}{2}n^2-\frac{29}{3}n
$;

$W(M_n)=5\sum\limits_{k=1}^{n-1}(n-k)(10(k-1)+9)+\frac{45}{2}n^2+\frac{9}{2}n=\frac{25}{3}n^3+20n^2-\frac{4}{3}n
$;

$W(P_n)=5\sum\limits_{k=1}^{n-1}(n-k)(15(k-1)+9)+\frac{45}{2}n^2+\frac{9}{2}n=\frac{25}{2}n^3+\frac{15}{2}n^2+7n
$.

\vspace{3mm}

{\bf Corollary 2.3}. The Wiener indices of the spiro orth-chain
$O_n$, the spiro meta-chain $M_n$ and the spiro para-chain $P_n$ are
$$W(O_n)=\frac{25}{6}n^3+\frac{65}{2}n^2-\frac{29}{3}n;$$
$$W(M_n)=\frac{25}{3}n^3+20n^2-\frac{4}{3}n;$$
$$W(P_n)=\frac{25}{2}n^3+\frac{15}{2}n^2+7n.$$

\vspace{3mm}

In the following, we consider the extremal problems of Wiener
indices among all spiro hexagonal chains with $n$ hexagons.

Let $G_n=H_0H_1\cdots H_{n-1}$ be a spiro hexagonal chain with $n$
hexagons, $c_{k}$ is the common cut-vertex of $H_{k-1}$ and $H_{k}$,
$1\leq k\leq n-1$. From Theorem 2.2 and $5k+9<10k+9<15k+9$, i.e.,
$f(o_k)<f(m_k)<f(p_k)$ for $k>1$, it is easily showed that $O_n$ is
the unique spiro hexagonal chain with the minimum Wiener index, and
the unique spiro hexagonal chain with the second minimal Wiener
index is the spiro hexagonal chain $G_n$ with the cut-vertex
sequence $(c_2, \cdots, c_{n-2}, c_{n-1})=(o_1, \cdots, o_{n-3},
m_{n-2})$. In order to find the third minimal value, we only need to
compare $2f(m_{n-3})+f(o_{n-2})$ with $2f(o_{n-3})+f(p_{n-2})$ from
Theorem 2.2. Since $2f(m_{n-3})+f(o_{n-2})<2f(o_{n-3})+f(p_{n-2})$,
the unique spiro hexagonal chain with the third minimal Wiener index
is the spiro hexagonal chain $G_n$ with $(c_2, \cdots, c_{n-3},
c_{n-2}, c_{n-1})=(o_1, \cdots, o_{n-4}, m_{n-3}, o_{n-2})$.

\vspace{3mm}

{\bf Theorem 2.4}. Among all spiro hexagonal chains with $n (n\geq
4)$ hexagons, (i) the unique spiro hexagonal chain with the minimum
Wiener index is $O_n$; (ii) the unique spiro hexagonal chain with
the second minimal Wiener index is the spiro hexagonal chain $G_n$
with the cut-vertex sequence $(c_2, \cdots, c_{n-1})=(o_1, \cdots,
o_{n-3}, m_{n-2})$; (iii) the unique spiro hexagonal chain with the
third minimal Wiener index is the spiro hexagonal chain $G_n$ with
$(c_2, \cdots, c_{n-1})=(o_1, \cdots, o_{n-4}, m_{n-3}, o_{n-2})$.

\vspace{3mm}

Analogously, the following results can be obtained. We omit their
proof and leave it for the reader.

\vspace{3mm}

{\bf Theorem 2.5}. Among all spiro hexagonal chains with $n (n\geq
4)$ hexagons, (i) the unique spiro hexagonal chain with the maximum
Wiener index is $P_n$; (ii) the unique spiro hexagonal chain with
the second maximal Wiener index is the spiro hexagonal chain $G_n$
with the cut-vertex sequence $(c_2, \cdots, c_{n-2}, c_{n-1})=(p_1,
\cdots, p_{n-3}, m_{n-2})$; (iii) the unique spiro hexagonal chain
with the third maximal Wiener index is the spiro hexagonal chain
$G_n$ with $(c_2, \cdots, c_{n-3}, c_{n-2}, c_{n-1})=(p_1, \cdots,
p_{n-4}, m_{n-3}, p_{n-2})$.

\section{Wiener index of polyphenyl hexagonal chains}

In this section, we will give a recurrence for computing the Wiener
indices of polyphenyl hexagonal chains, and then derive a formula
for computing the Wiener indices of polyphenyl hexagonal chains.

Let $G$ be a cactus graph in which each block is either an edge or a
hexagon. $G$ is called a polyphenyl hexagonal chain if each hexagon
of $G$ has at most two cut-vertices, and each cut-vertex is shared
by exactly one hexagon and one cut-edge. The number of hexagons in
$G$ is called the length of $G$. An example of a polyphenyl
hexagonal chain is shown in Figure 2(ii).

Obviously, a polyphenyl hexagonal chain of length $n$ has $6n$
vertices and $7n-1$ edges. Furthermore, any polyphenyl hexagonal
chain of length greater than one has exactly two hexagons with only
one cut-vertex. Such hexagons are called terminal hexagons. Any
remaining hexagons are called internal hexagons.

Note that any polyphenyl hexagonal chain
$\overline{G}_n=\overline{H}_0\overline{H}_1\cdots
\overline{H}_{n-1}$ of length $n (n\geq 2)$ can be obtained from the
polyphenyl hexagonal chain
$\overline{G}_{n-1}=\overline{H}_0\overline{H}_1\cdots
\overline{H}_{n-2}$ of length $n-1$ by a cut-edge linking a vertex
$c_{n-1}$ in the hexagon $\overline{H}_{n-1}$ to a non cut-vertex
$u$ in the terminal hexagon $\overline{H}_{n-2}$ of
$\overline{G}_{n-1}$, where $u$ is said to be the tail of
$\overline{H}_{n-1}$, denoted by $t_{n-1}$. A vertex $v$ of
$\overline{H}_{n-1}$ is said to be ortho-, meta- and para-vertex if
the distance between $v$ and $c_{n-1}$ is 1, 2 and 3, denoted by
$o_{n-1}$, $m_{n-1}$ and $p_{n-1}$, respectively. Examples of tail,
ortho-, meta-, and para-vertices are shown in Figure 5.

A polyphenyl hexagonal chain $\overline{G}_n$ is a polyphenyl
ortho-chain if $t_k=o_{k-1}$ for $2\leq k\leq n-1$. The polyphenyl
meta-chain and polyphenyl para-chain are defined in a completely
analogous manner. The polyphenyl ortho-, meta- and para-chain of
length $n$ is denoted by $\overline{O}_n$, $\overline{M}_n$ and
$\overline{P}_n$, respectively. Examples of polyphenyl ortho-,
meta-, and para-chains are shown in Figure 6.

Let $\overline{G}_n=\overline{H}_0\overline{H}_1\cdots
\overline{H}_{n-1}$ be a polyphenyl hexagonal chain with $n$
hexagons. $\overline{G}_{n-1}=\overline{H}_0\overline{H}_1\cdots
\overline{H}_{n-2}$  and $t_{n-1}, o_{n-1}, m_{n-1}, p_{n-1}$ are
the tail-, ortho-, meta-, and para-vertex in $\overline{H}_{n-1}$,
respectively. Then

$$
\begin{array}{ll}
W(\overline{G}_n, o_{n-1})&=\sum\limits_{v\in \overline{G}_n}d(v,
o_{n-1})=\sum\limits_{v\in \overline{G}_{n-1}}(d(v, t_{n-1})+2)+9 \\
&=W(\overline{G}_{n-1}, t_{n-1})+12(n-1)+9;\\
\end{array}
$$

$$
\begin{array}{ll}
W(\overline{G}_n, m_{n-1})&=\sum\limits_{v\in \overline{G}_n}d(v,
m_{n-1})=\sum\limits_{v\in \overline{G}_{n-1}}(d(v, t_{n-1})+3)+9\\
&=W(\overline{G}_{n-1}, t_{n-1})+18(n-1)+9;\\
\end{array}
$$

$$
\begin{array}{ll}
W(\overline{G}_n, p_{n-1})&=\sum\limits_{v\in \overline{G}_n}d(v,
p_{n-1})=\sum\limits_{v\in \overline{G}_{n-1}}(d(v, t_{n-1})+4)+9\\
&=W(\overline{G}_{n-1}, t_{n-1})+24(n-1)+9.\\
\end{array}
$$

So, we have
\begin{equation}
W(\overline{G}_n, t_n)=W(\overline{G}_{n-1}, t_{n-1})+g(t_n)
\end{equation}
where
$$
g(t_n)=\left\{
\begin{array}{ll}
12(n-1)+9, & t_n \,\,\, is \,\,\, o_{n-1};\\
18(n-1)+9, & t_n \,\,\, is \,\,\, m_{n-1};\\
24(n-1)+9, & t_n \,\,\, is \,\,\, p_{n-1}.\\
\end{array}
\right.
$$

\vspace{3mm}

Let $v_1, v_2, v_3, v_4, v_5, v_6$ be the vertices of
$\overline{H}_{n-1}$ different from the tail $t_{n-1}$. By the
definition of Wiener index,
$$
\begin{array}{ll}
W(\overline{G}_n)&=W(\overline{G}_{n-1})+\sum\limits_{i=1}^6\sum\limits_{v\in
\overline{G}_{n-1}}d(v,
v_i)+\sum\limits_{1\leq i<j\leq 6}d(v_i, v_j)\\
&=W(\overline{G}_{n-1})+\sum\limits_{v\in \overline{G}_{n-1}}(6d(v,
t_{n-1})+15)+27\\
&=W(\overline{G}_{n-1})+6W(\overline{G}_{n-1}, t_{n-1})+90n-63,\\
\end{array}
$$
i.e.,
\begin{equation}
W(\overline{G}_n)=W(\overline{G}_{n-1})+6W(\overline{G}_{n-1},
t_{n-1})+90n-63.
\end{equation}

\vspace{3mm}

Combining equations (6) and (7), we can get the following recurrence
for computing the Wiener indices of polyphenyl hexagonal chains.

{\bf Theorem 3.1}. Let
$\overline{G}_n=\overline{H}_0\overline{H}_1\cdots
\overline{H}_{n-1}$ be a polyphenyl hexagonal chain with $n$
hexagons, $t_{n-1}$ the tail of $\overline{H}_{n-1}$. Then
$$
\left\{
\begin{array}{l}
W(\overline{G}_n)=W(\overline{G}_{n-1})+6W(\overline{G}_{n-1}, t_{n-1})+90n-63, \\
W(\overline{G}_n, t_n)=W(\overline{G}_{n-1}, t_{n-1})+g(t_n),\\
W(\overline{G}_1)=W(\overline{H}_0)=27, \\
W(\overline{G}_1, t_1)=g(t_1)=9,\\
\end{array}
\right.
$$ and
$$
g(t_n)=\left\{
\begin{array}{ll}
12(n-1)+9, & t_n \,\,\, is \,\,\, o_{n-1};\\
18(n-1)+9, & t_n \,\,\, is \,\,\, m_{n-1};\\
24(n-1)+9, & t_n \,\,\, is \,\,\, p_{n-1}.\\
\end{array}
\right.
$$

\vspace{3mm}

\vspace{3mm}

Using the recurrence above, we have
$$
\begin{array}{ll}
W(\overline{G}_n)&=W(\overline{G}_{n-1})+6W(\overline{G}_{n-1}, t_{n-1})+90n-63\\
&=W(\overline{G}_1)+6\sum\limits_{k=1}^{n-1}W(\overline{G}_{n-k},
t_{n-k})+90\sum\limits_{k=1}^{n-1}(n-k+1)-63(n-1)\\
&=27+6\sum\limits_{k=1}^{n-1}W(\overline{G}_{k},
t_{k})+45(n+2)(n-1)-63(n-1)\\
&=6\sum\limits_{k=1}^{n-1}W(\overline{G}_{k}, t_{k})+45n^2-18n,\\
\end{array}
$$
i.e.,
\begin{equation}
W(\overline{G}_n)=6\sum\limits_{k=1}^{n-1}W(\overline{G}_{k},
t_{k})+45n^2-18n.
\end{equation}

\vspace{3mm}

And,

$$
\begin{array}{ll}
W(\overline{G}_k, t_k)&=W(\overline{G}_{k-1}, t_{k-1})+g(t_k)\\
&=W(\overline{G}_1, t_1)+g(t_2)+\cdots+g(t_k)\\
&=g(t_1)+g(t_2)+\cdots+g(t_k),
\end{array}
$$
i.e.,
\begin{equation}
W(\overline{G}_k, t_k)=\sum\limits_{i=1}^{k}g(t_{i}).
\end{equation}

Combining equations (8) and (9), we can obtain the following formula
for computing the Wiener indices of polyphenyl hexagonal chains.

\vspace{3mm}

{\bf Theorem 3.2}. Let
$\overline{G}_n=\overline{H}_0\overline{H}_1\cdots
\overline{H}_{n-1}$ be a polyphenyl hexagonal chain with $n$
hexagons, $t_{k}$ the tail of $\overline{H}_{k}$, $1\leq k\leq n-1$.
Then
$$
\begin{array}{ll}
W(\overline{G}_n)&=6\sum\limits_{k=1}^{n-1}\sum\limits_{i=1}^kg(t_{i})+45n^2-18n\\
&=6\sum\limits_{k=1}^{n-1}(n-k)g(t_k)+45n^2-18n,\\
\end{array}
$$
where
\begin{equation}
g(t_k)=\left\{
\begin{array}{ll}
12(k-1)+9, & t_k\,\,\, is \,\,\, o_{k-1};\\
18(k-1)+9, & t_k \,\,\, is \,\,\, m_{k-1};\\
24(k-1)+9, & t_k \,\,\, is \,\,\, p_{k-1}.\\
\end{array}
\right. \end{equation} and $g(t_1)=9$.

\vspace{3mm}

For the polyphenyl orth-chain $\overline{O}_n$, the polyphenyl
meta-chain $\overline{M}_n$ and the polyphenyl para-chain
$\overline{L}_n$, $t_k$ is $o_{k-1}, m_{k-1}$ and $p_{k-1}$,
respectively. So, we have

$W(\overline{O}_n)=6\sum\limits_{k=1}^{n-1}(n-k)(12(k-1)+9)+45n^2-18n=12n^3+36n^2-21n
$;

$W(\overline{M}_n)=6\sum\limits_{k=1}^{n-1}(n-k)(15(k-1)+9)+45n^2-18n=18n^3+18n^2-9n
$;

$W(\overline{P}_n)=6\sum\limits_{k=1}^{n-1}(n-k)(24(k-1)+9)+45n^2-18n=24n^3+3n
$.

\vspace{3mm}

{\bf Corollary 3.3}. The Wiener indices of the polyphenyl orth-chain
$\overline{O}_n$, the polyphenyl meta-chain $\overline{M}_n$ and the
polyphenyl para-chain $\overline{L}_n$ are
$$W(\overline{O}_n)=12n^3+36n^2-21n;$$
$$W(\overline{M}_n)=18n^3+18n^2-9n;$$
$$W(\overline{P}_n)=24n^3+3n.$$

\section{A relation between $W(G_n)$ and $W(\overline{G}_n)$}

An exact relation between the Wiener indices of a phenylene and its
hexagonal squeeze was established by Pavlovi\'{c} and Gutman [13].

To every polyhenyl hexagonal chain, it is possible to associate a
spiro hexagonal chain, obtained so that the cut edges of the
polyphenyl hexagonal chain are squeezed off. This spiro hexagonal
chain is named the hexagonal squeeze of the respective polyphenyl
hexagonal chain. Clearly, each polyhenyl hexagonal chain determines
a unique hexagonal squeeze and vice versa, and these two systems
have an equal number of hexagons. For example, the spiro hexagonal
chain in Figure 2(i) is the hexagonal squeeze of the polyphenyl
hexagonal chain in Figure 2(ii). Here, we also give a relation
between the Wiener indices of a polyphenyl hexagonal chain and its
hexagonal squeeze.

{\bf Theorem 4.1}. Let
$\overline{G}_n=\overline{H}_0\overline{H}_1\cdots
\overline{H}_{n-1}$ be a polyphenyl hexagonal chain with $n$
hexagons, $G_n=H_0H_1\cdots H_{n-1}$ its hexagonal squeeze. The
Wiener indices of $\overline{G}_n$ and $G_n$ are related as
\begin{equation}
25W(\overline{G}_n) =36W(G_n)+150n^3-270n^2-177n.
\end{equation}

{\bf Proof}. From equations (5) and (10), we have
$$\frac{g(t_k)-9}{6}=\frac{f(c_k)-9}{5}+(k-1),$$
i.e.,$$5g(t_k)=6f(c_k)+30k-39.$$ So,
$$5\sum\limits_{k=1}^{n-1}(n-k)g(t_k)
=6\sum\limits_{k=1}^{n-1}(n-k)f(c_k)+\sum\limits_{k=1}^{n-1}(n-k)(30k-39).$$

By Theorems 2.2 and 3.2,
$$\frac{5}{6}W(\overline{G}_n)-\frac{5}{6}(45n^2-18n)
=\frac{6}{5}W(G_n)-\frac{6}{10}(45n^2+9n)+\sum\limits_{k=1}^{n-1}(n-k)(30k-39),$$
i.e.,
$$\frac{5}{6}W(\overline{G}_n)
=\frac{6}{5}W(G_n)+5n^3-9n^2-\frac{59}{10}n$$ and
$$25W(\overline{G}_n)
=36W(G_n)+150n^3-270n^2-177n.$$

\vspace{3mm}

From Theorem 4.1, we can obtain the following result on the extremal
problems of polyphenyl hexagonal chains with respect to the Wiener
index.

\vspace{3mm}

{\bf Theorem 4.2}. (i) Among all polyphenyl hexagonal chains with $n
(n\geq 4)$ hexagons, $\overline{G}_n$ has the minimum, the second
and the third minimal Wiener index if and only if its hexagonal
squeeze $G_n$ has the minimum, the second and the third minimal
Wiener index among all spiro hexagonal chains with $n$ hexagons.

(ii) Among all polyphenyl hexagonal chains with $n (n\geq 4)$
hexagons, $\overline{G}_n$ has the maximum, the second and the third
maximal Wiener index if and only if its hexagonal squeeze $G_n$ has
the maximum, the second and the third maximal Wiener index among all
spiro hexagonal chains with $n$ hexagons.

\section{The average value of the Wiener index}

If $\mathcal{G}_n$ is the set of all spiro hexagonal chains with $n$
hexagons, then the average value of the Wiener indices with respect
to $\mathcal{G}_n$ is
$$W_{avr}(\mathcal{G}_n)=\frac{1}{|\mathcal{G}_n|}\sum\limits_{G\in
\mathcal{G}_n}W(G).$$

Let $G_n=H_0H_1\cdots H_{n-1}$ be a spiro hexagonal chain of length
$n$. $c_k$ is the common cut-vertex of $H_{k-1}$ and $H_k$,
$k=1,2,\cdots,n-1$. $H_k$ is called ortho-hexagon, meta-hexagon, or
para-hexagon in [9] if the distance between its cut-vertices
$c_{k-1}$ and $c_k$ is 1, 2, and 3, respectively. In an obvious way,
each spiro hexagonal chain of length $n$ can be assigned a word of
length $n-2$ over the alphabet $\{O, M, P\}$. Such a word is called
the code of the chain. For example, the code of the chain in Figure
3(i) is PMMMO. The correspondence is not necessarily bijective: the
same chain is also described by the code OMMMP, i.e., the same code
read backwards. It is easy to see that a palindromic code uniquely
defines a chain, while exactly two non-palindromic codes correspond
to the same chain. From here, it was concluded in [9] that the
number of all possible spiro hexagonal chains of length $n$ is equal
to the number obtained by adding half the number of non-palindromic
codes of length $n-2$ to the number of palindromic codes of the same
length. Since the number of palindromic codes of length $n$ is equal
to $3^{\lfloor\frac{n+1}{2}\rfloor}$ and the total number of codes
is equal to $3^n$, we have the following result.

\vspace{3mm}

{\bf Lemma 5.1}([9]). The number of different spiro hexagonal chains
of length $n$ is
$$|\mathcal{G}_n|=\frac{1}{2}(3^{n-2}+3^{\lfloor\frac{n-1}{2}\rfloor}).$$

\vspace{3mm}

Since $f(c_1)=9$, from Theorem 2.2, the Wiener index of a spiro
hexagonal chain $G_n$ of length $n$ can be reduced to
$$
\begin{array}{ll}
W(G_n)&=5\sum\limits_{k=2}^{n-1}(n-k)f(c_k)+\frac{9}{2}(5n^2+11n-10).\\
\end{array}
$$

Let $G_n=H_0H_1\cdots H_{n-1}$ be a spiro hexagonal chain of length
$n$. $x_1x_2\cdots x_{n-2}$ is its code, where $x_i\in\{O, M, P\}$.
$(c_2, c_3, \cdots, c_{n-1})$ is its cut-vertex sequence.
Then
$$\left\{
\begin{array}{lll}
x_i=O & \mbox{if and only if} & c_{i+1}=o_i; \\
x_i=M & \mbox{if and only if} & c_{i+1}=m_i; \\
x_i=P & \mbox{if and only if} & c_{i+1}=p_i. \\
\end{array}
\right.
$$

Note that each of $o_{k-1}, m_{k-1}, p_{k-1}$ can be taken $3^{n-3}$
times by $c_k$ when $(c_2, \cdots, c_{n-1})$ is taken over all the
sequences of length $n-2$, and each of $o_{k-1}, m_{k-1}, p_{k-1}$
can be taken by $3^{\lfloor\frac{n-3}{2}\rfloor}$ times $c_k$ when
$(c_2, \cdots, c_{n-1})$ is taken over all the palindromic sequences
of length $n-2$. We have the following result
\begin{center}
$\sum\limits_{G_n\in\mathcal{G}_n}W(G_n)=\frac{1}{2}(\sum_1
W(G_n)+\sum_2 W(G_n))$ \end{center} where $\sum_1$ is taken over all
$G_n$ whose cut-vertex sequences $(c_2, \cdots, c_{n-1})$ are the
sequences of length $n-2$, and $\sum_2$ is taken over all $G_n$
whose cut-vertex sequences $(c_2, \cdots, c_{n-1})$ are the
palindromic sequences of length $n-2$. So,
$$
\begin{array}{lll}
\sum_1W(G_n)&=&
3^{n-3}\times 5\sum\limits_{k=2}^{n-1}(n-k)[f(o_{k-1})+f(m_{k-1})+f(p_{k-1})]\\
&&+\frac{9}{2}(5n^2+11n-10)\times 3^{n-2}\\
&=&3^{n-3}\times 5\sum\limits_{k=2}^{n-1}(n-k)[30(k-1)+27]+\frac{3^n}{2}(5n^2+11n-10)\\
&=&3^{n-3}\times 5\times(5n^3-\frac{3}{2}n^2-\frac{61}{2}n+27)+\frac{3^n}{2}(5n^2+11n-10)\\
&=&3^{n-3}(25n^3+60n^2-4n)\\
\end{array}
$$
$$
\begin{array}{lll}
\sum_2W(G_n)&=&
3^{\lfloor\frac{n-3}{2}\rfloor}\times 5\sum\limits_{k=2}^{n-1}(n-k)[f(o_{k-1})+f(m_{k-1})+f(p_{k-1})]\\
&&+\frac{9}{2}(5n^2+11n-10)\times 3^{\lfloor\frac{n-1}{2}\rfloor}\\
&=&3^{\lfloor\frac{n-3}{2}\rfloor}(25n^3+60n^2-4n)\\
\end{array}
$$
and
$$\sum\limits_{G_n\in\mathcal{G}_n}W(G_n)=\frac{1}{2}(3^{n-3}+3^{\lfloor\frac{n-3}{2}\rfloor})(25n^3+60n^2-4n).$$

By Lemma 5.1, we can get the average value of the Wiener indices
with respect to $\mathcal{G}_n$.

\vspace{3mm}

{\bf Theorem 5.2}. The average value of the Wiener indices with
respect to $\mathcal{G}_n$ is
$$W_{avr}(\mathcal{G}_n)=\frac{1}{3}(25n^3+60n^2-4n).$$

Note that the average value of the Wiener indices with respect to
$\{O_n, M_n, P_n\}$ is
$$\frac{W(O_n)+W(M_n)+W(P_n)}{3}=\frac{1}{3}(25n^3+60n^2-4n)$$
from Corollary 2.3. The interesting result shows that the average
value of the Wiener indices with respect to $\mathcal{G}_n$ is
exactly the average value of the Wiener indices with respect to
$\{O_n, M_n, P_n\}$, and is just the Wiener index $W(M_n)$ of the
spiro meta-chain $M_n$.

\vspace{3mm}

Similarly, if $\mathcal{\overline{G}}_n$ is the set of all
polyphenyl hexagonal chains with $n$ hexagons, then the average
value of the Wiener indices with respect to
$\mathcal{\overline{G}}_n$ is
$$W_{avr}(\mathcal{\overline{G}}_n)=\frac{1}{|\mathcal{\overline{G}}_n|}\sum\limits_{G\in
\mathcal{\overline{G}}_n}W(G).$$

Using the hexagonal squeeze, there is a bijection between
$\mathcal{\overline{G}}_n$ and $\mathcal{G}_n$. So, we have
$$|\mathcal{\overline{G}}_n|=|\mathcal{G}_n|=\frac{1}{2}(3^{n-2}+3^{\lfloor\frac{n-1}{2}\rfloor})$$
and
\begin{center}
$\sum\limits_{\overline{G}_n\in\mathcal{\overline{G}}_n}W(\overline{G}_n)=\frac{1}{2}(\sum_1
W(\overline{G}_n)+\sum_2 W(\overline{G}_n))$ \end{center} where
$\sum_1$ is taken over all $\overline{G}_n$ whose cut-vertex
sequences $(c_2, \cdots, c_{n-1})$ of their hexagonal squeezes are
the sequences of length $n-2$, and $\sum_2$ is taken over all
$\overline{G}_n$ whose cut-vertex sequences $(c_2, \cdots, c_{n-1})$
of their hexagonal squeezes are the palindromic sequences of length
$n-2$.

Since $g(t_1)=9$, from Theorem 3.2, the Wiener index of a polyphenyl
hexagonal chain $\overline{G}_n$ of length $n$ can be reduced to
$$
\begin{array}{ll}
W(\overline{G}_n)&=6\sum\limits_{k=2}^{n-1}(n-k)g(t_k)+(45n^2+36n-54).\\
\end{array}
$$

$$
\begin{array}{lll}
\sum_1W(\overline{G}_n)&=&
3^{n-3}\times 6\sum\limits_{k=2}^{n-1}(n-k)[g(o_{k-1})+g(m_{k-1})+g(p_{k-1})]\\
&&+3^{n-2}(45n^2+36n-54)\\
&=&3^{n-3}\times 6\sum\limits_{k=2}^{n-1}(n-k)[54(k-1)+27]+3^{n-2}(45n^2+36n-54)\\
&=&3^{n-3}(54n^3-81n^2-135n+162)+3^{n-2}(45n^2+36n-54)\\
&=&3^{n-3}(54n^3+54n^2-27n)\\
\end{array}
$$
$$
\begin{array}{lll}
\sum_2W(\overline{G}_n)&=&
3^{\lfloor\frac{n-3}{2}\rfloor}\times 6\sum\limits_{k=2}^{n-1}(n-k)[g(o_{k-1})+g(m_{k-1})+g(p_{k-1})]\\
&&+3^{\lfloor\frac{n-1}{2}\rfloor}(45n^2+36n-54)\\
&=&3^{\lfloor\frac{n-3}{2}\rfloor}(54n^3+54n^2-27n)\\
\end{array}
$$
and
$$\sum\limits_{\overline{G}_n\in\mathcal{G}_n}W(G_n)=\frac{1}{2}(3^{n-3}
+3^{\lfloor\frac{n-3}{2}\rfloor})(54n^3+54n^2-27n).$$

Hence, we can get the average value of the Wiener indices with
respect to $\mathcal{\overline{G}}_n$.

\vspace{3mm}

{\bf Theorem 5.3}. The average value of the Wiener indices with
respect to $\mathcal{\overline{G}}_n$ is
$$W_{avr}(\mathcal{\overline{G}}_n)=18n^3+18n^2-9n.$$

Note that the average value of the Wiener indices with respect to
$\{\overline{O}_n, \overline{M}_n, \overline{P}_n\}$ is
$$\frac{W(\overline{O}_n)+W(\overline{M}_n)+W(\overline{P}_n)}{3}=18n^3+18n^2-9n$$
from Corollary 3.3. This also shows that the average value of the
Wiener indices with respect to $\mathcal{\overline{G}}_n$ is exactly
that to $\{\overline{O}_n, \overline{M}_n, \overline{P_n}\}$, and is
just the Wiener index of the polyphenyl meta-chain $\overline{M}_n$.

\vspace{3mm}

Finally, by Theorems 5.2 and 5.3, we have
$$25W_{avr}(\overline{\mathcal{G}}_n)
=36W_{avr}(\mathcal{G}_n)+150n^3-270n^2-177n.$$ This relation is
identical with the equation (11) in Theorem 4.1.

\includegraphics{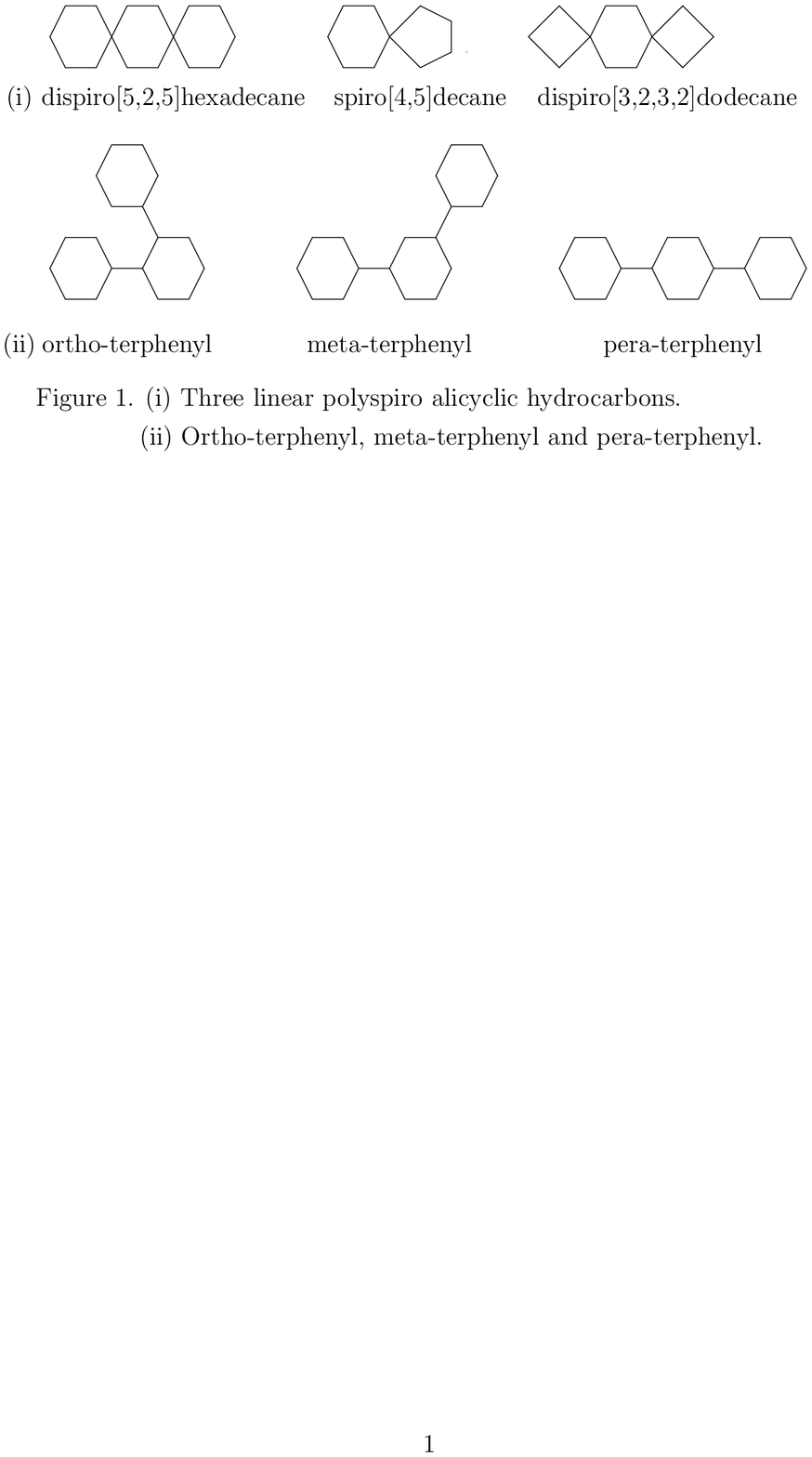}

\includegraphics{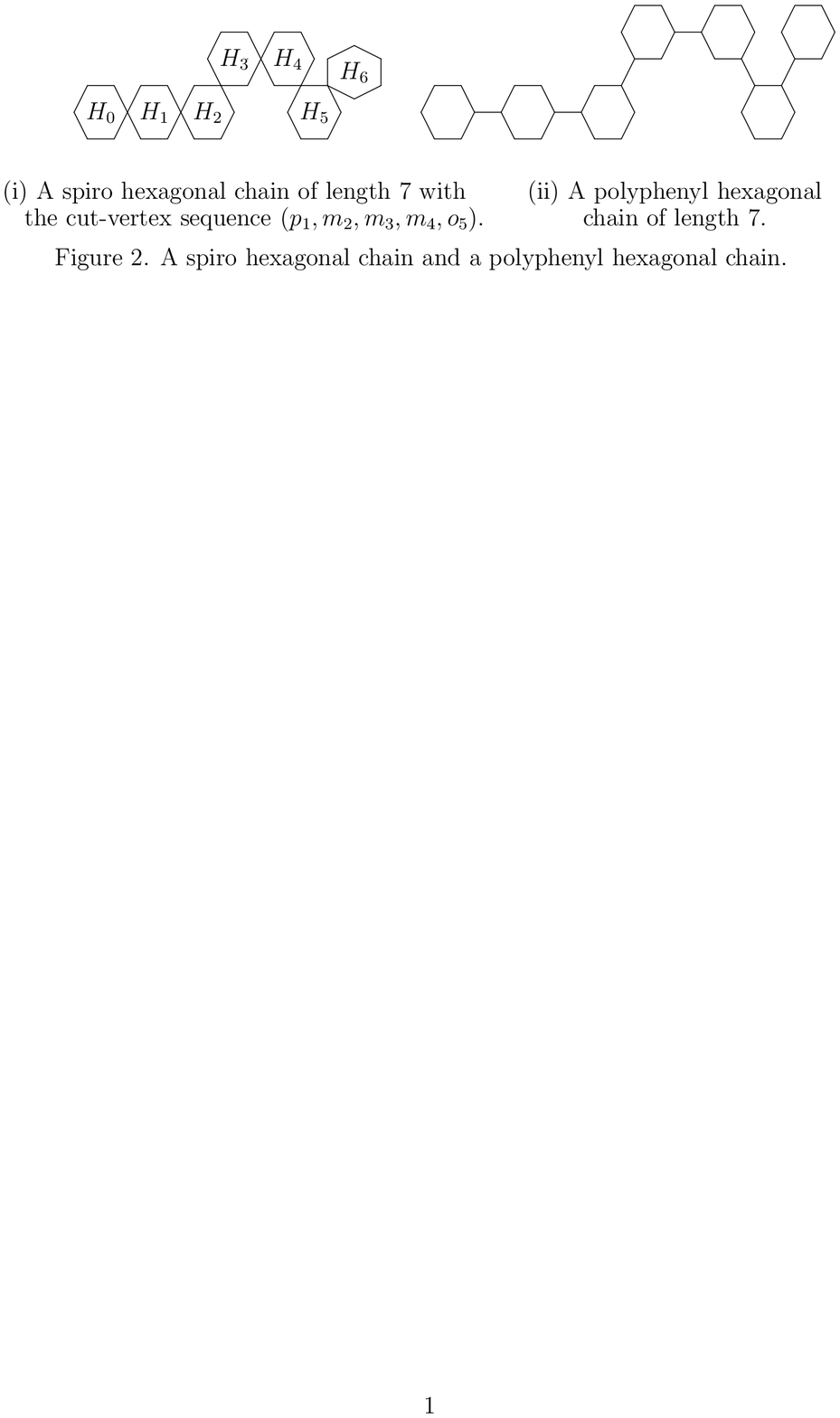}

\includegraphics{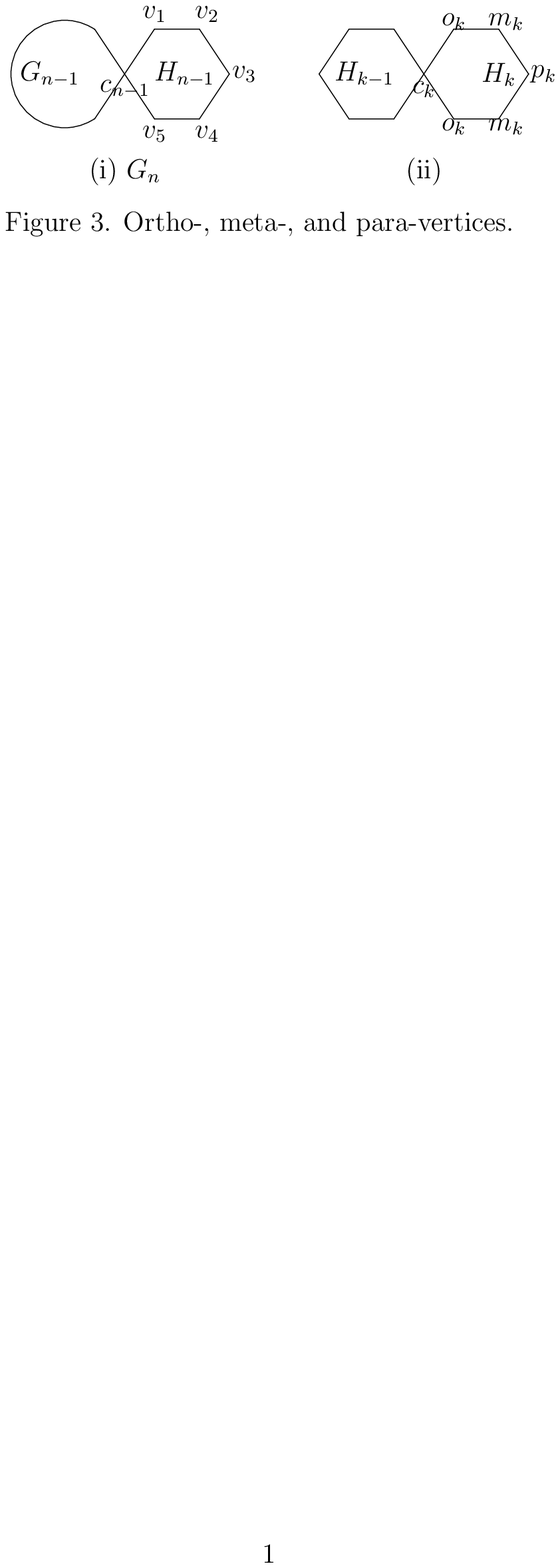}

\includegraphics{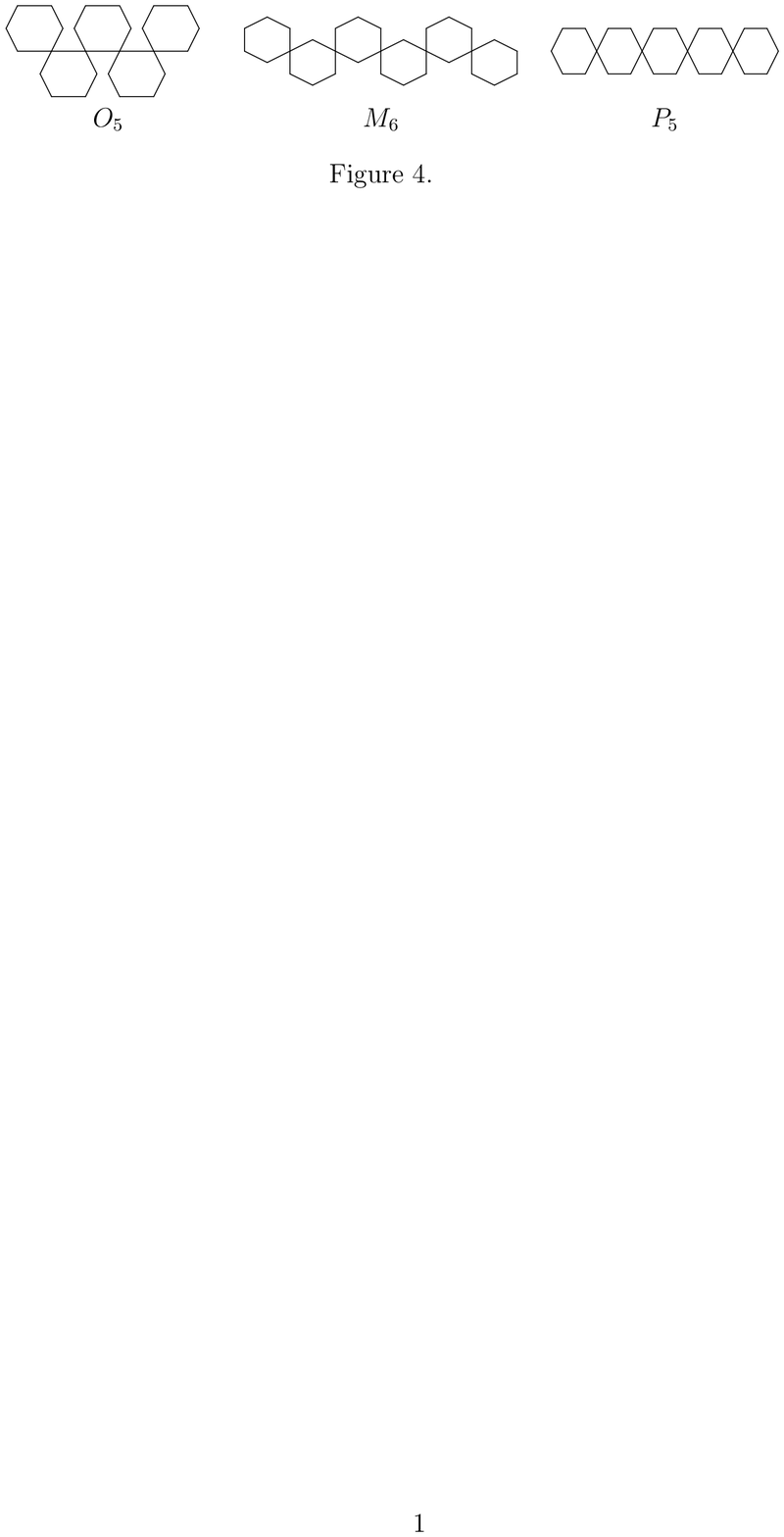}

\includegraphics{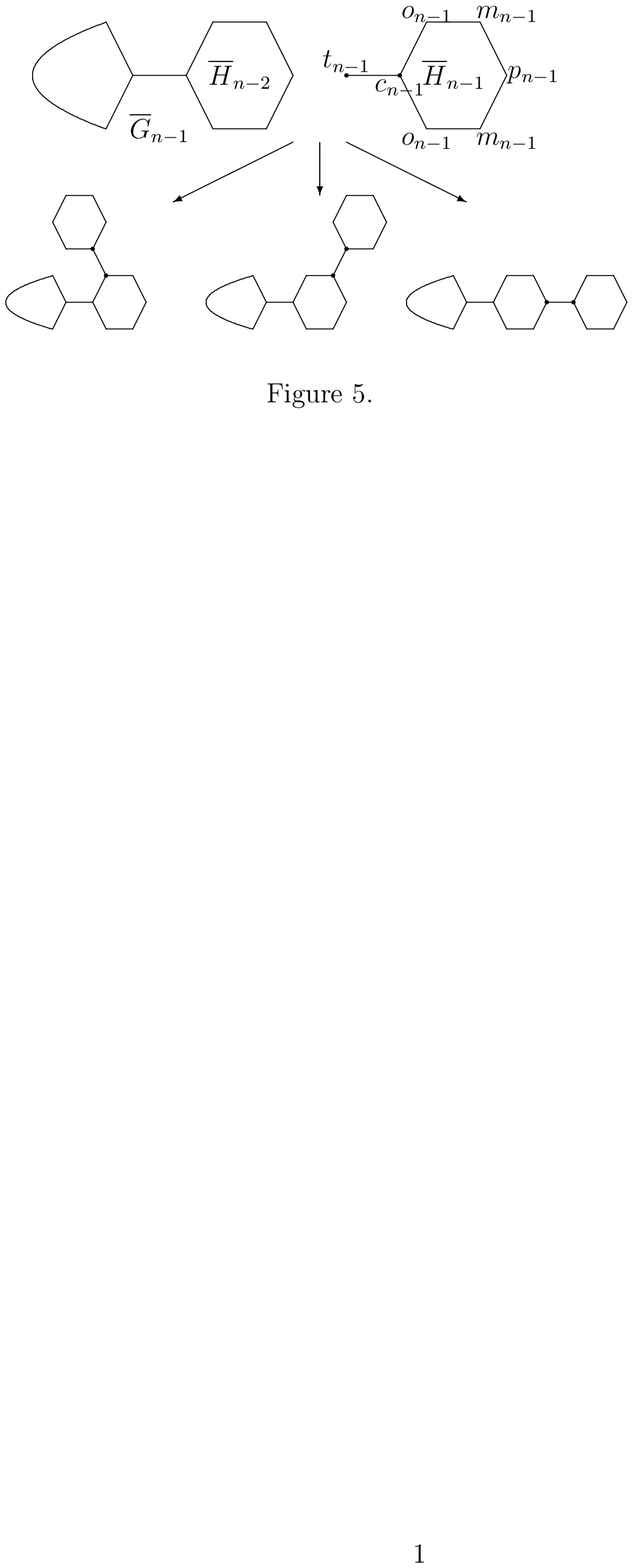}

\includegraphics{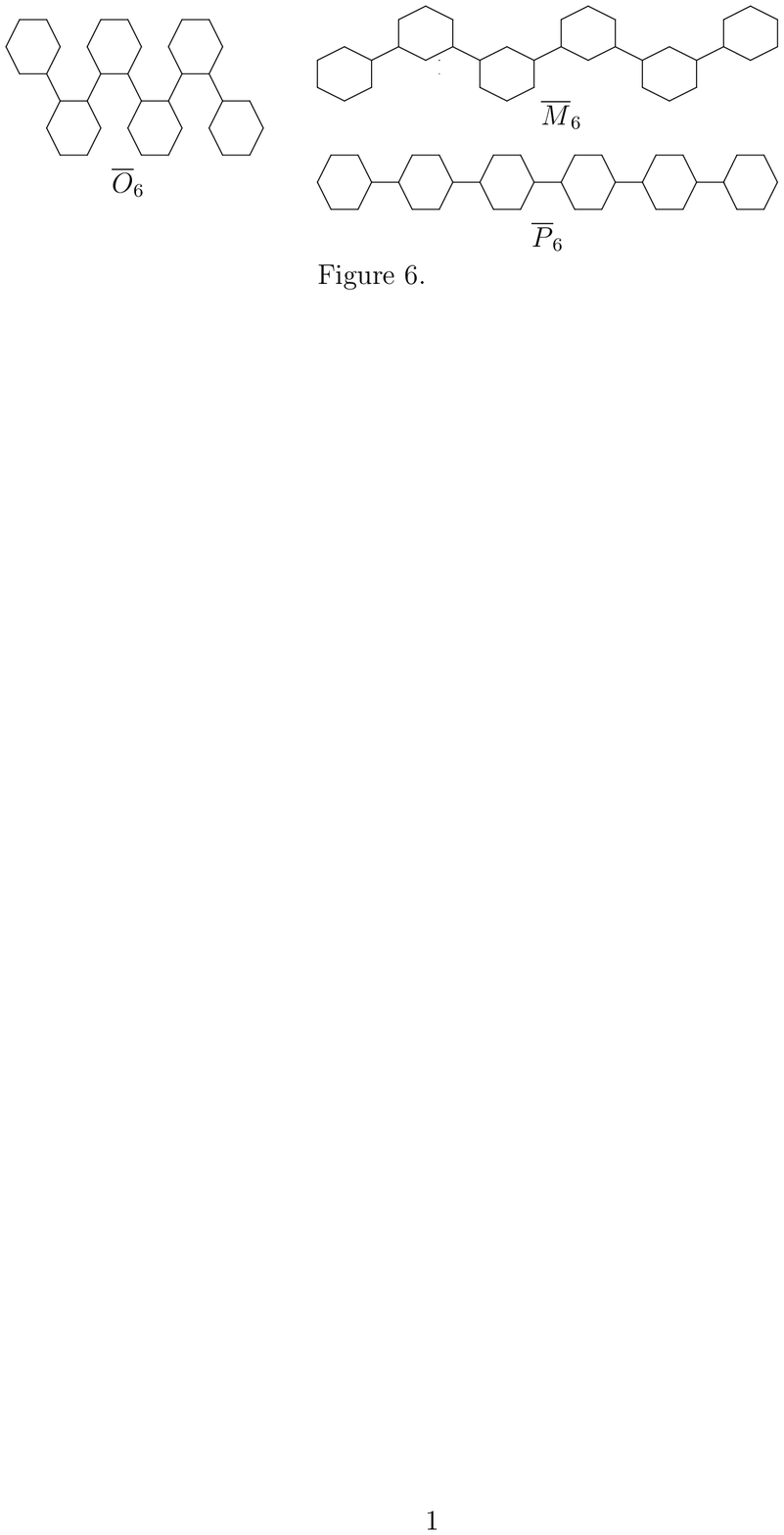}

\end{document}